\documentclass[12pt]{article}
\usepackage[cp1251]{inputenc}
\usepackage[T2A]{fontenc}
\usepackage[english]{babel}

\def\eps{\varepsilon}

\newcounter{num}[section]

\newcommand{\Th}{\refstepcounter{num}
{\bf Theorem \arabic{section}.\arabic{num} }}

\newcommand{\Lemma}{\refstepcounter{num}
{\bf Lemma \arabic{section}.\arabic{num} }}

\newcommand{\Cor}{\refstepcounter{num}
{\bf Corollary \arabic{section}.\arabic{num} }}

\newcommand{\Def}{\refstepcounter{num}
{\it Definition \arabic{section}.\arabic{num} }}

\newcommand{\Proof}{{\bf Proof. }}

\def\_phi{\varphi}
\def\d{\delta}

\author{Shkredov I.D.}
\title{On multiple recurrence.}
\date{}
\begin{document}
\maketitle

\refstepcounter{section}
{\bf \arabic{section}. Introduction.}

 Let $N$ be a natural number and
$$
  a_k(N) = \frac{1}{N} \max \{ |A| ~:~ A \subseteq [1,N],
$$
$$
  A \mbox{ --- does not contain an arithmetic progression of length } k
   \},
$$
where $|A|$ denotes the cardinality of a set $A$.
    In \cite{EaT} P. Erdos and P. Turan realised that it ought
to be possible to find arithmetic progression of  length $k$
in any set with positive density.
   In other words they conjectured that for any $k\ge 3$
\begin{equation}
  a_k(N) \to 0, \mbox{ as } N\to \infty
\label{Sz_lab}
\end{equation}

  In case $k=3$ conjecture (\ref{Sz_lab}) was proved by
K.F. Roth in \cite{Rt}.
In his paper Roth used the Hardy -- Littlewood method
to prove the inequality
$$
  a_3(N) \ll \frac{1}{\log \log N}.
$$
  At this moment the best result about
  a lower bound  for $a_3(N)$
  belongs to J. Bourgain \cite{Bu}.
He proved that
\begin{equation}
 a_3(N) \ll \sqrt{ \frac{\log \log N}{ \log N} }.
\end{equation}

  For an arbitrary $k$ conjecture (\ref{Sz_lab}) was proved by
E. Szemeredi \cite{Sz} in 1975.

  The second proof of Szemeredi's theorem was given by
 H. Furstenberg in
\cite{Fu}, using ergodic theory.
Furstenberg showed that
Szemeredi's theorem is equivalent to the multiple recurrence of
almost every point in an arbitrary dynamical system.
Here we formulate this theorem in the case of metric spaces :
\\
\Th {\it Let $X$ be a metric space with metric
$d(\cdot,\cdot)$ and Borel sigma--algebra of measurable sets
$\Phi$. Let $T$ be
a measurable map of {\it X} into itself preserving the measure $\mu$
and let $k\ge 3$.
Then
$$
  \liminf_{n\to \infty} \max
  \{ d(T^{n}x,x), d(T^{2n}x,x), \dots, d(T^{(k-1)n}x,x) \} = 0.
$$
for almost all $x\in X$.
}\label{Furst}

Actually, H. Furstenberg obtained  more general result.

\Th {\it Let $X$ be a metric space with metric
$d(\cdot,\cdot)$ and Borel sigma--algebra of measurable sets
$\Phi$. Let $l\in {\bf N}$ and $T_1,\dots, T_l$ be commutative
measurable maps of {\it X} into itself preserving the measure $\mu$.
Then
$$
  \liminf_{n\to \infty} \max
  \{ d(T_1^{n}x,x), d(T_2^{2n}x,x), \dots, d(T_l^{n}x,x) \} = 0.
$$
for almost all $x\in X$.
}\label{Furst1}

  Unfortunately, Szemeredi's methods give very weak upper bound
for $a_k(N)$.
Furstenberg's proof gives no bound.
  Only in  2001
W.T. Gowers \cite{Gow_m}
obtained a quantitative result about the speed of tending to zero
of
$a_k(N)$ with $k\ge 4$.
He proved the following theorem.
\\
\Th {\it Let $\delta >0$, $k\ge 4$ and $N\ge \exp \exp (C
\delta^{-K})$, where
$C,K > 0$  is absolute constants.
Let $A\subseteq \{ 1,2,\dots ,N \} $ be a set
of cardinality at least
$\delta N$. Then $A$ contains an arithmetic progression
of length  $k$.}
\\
In other words, W.T. Gowers proved that for any $k\ge 4$, we have
$a_k(N) \ll 1/ (\log \log N)^{c_k}$, where constant $c_k$ depends
on $k$ only.

  Consider the following problem.
  Let us consider the two--dimensional lattice $[1,N]^2$ with basis  $\{(1,0)$, $(0,1)\}$.
Define
$$
  L(N) = \frac{1}{N^2} \max \{~ |A| ~:~ A\subseteq [1,N]^2 ~\mbox{ and }~
$$
$$
  A \mbox { --- does not contain any triple } \{ (k,m),~ (k+d,m),~ (k,m+d),~ d>0 \}
$$
\begin{equation}\label{tri}
    \mbox { with positive } d \}.
\end{equation}
 A triple from (\ref{tri}) will be called a
 {\it "corner".}
In papers \cite{Sz2,Fu}
shown that  $L(N)$ tends to   $0$ as $N$
tends to infinity.
W.T. Gowers (see \cite{Gow_m}) set a question
about the speed of convergence to $0$ of $L(N)$.
\\
In \cite{Vu} V. Vu proposed the following solution.
Let us define
$\log_{*}N$ as the largest integer $k$ such that
$\log_{[l]} N \ge 2$, where
$\log_{[1]} N = \log N$ and for $l\ge 2~$
$\log_{[l]} = \log( \log_{[l-1]} N)$.
V. Vu proved that
$$
  L(N) \le \frac{100}{\log_{*}^{1/4} N}
$$
The main result of \cite{Tri} is
\\
\Th {\it Let $\delta>0$, $N\ge \exp \exp \exp ( \d^{-c} )$,
where  $c>0$
is an
absolute constant.
Let
$A\subseteq \{1,\dots, N\}^2$
be a set
of cardinality at least $\delta N^2$.
Then $A$ contains a triple $(k,m), (k+d,m), (k,m+d)$, where
$d>0$.} \label{main_th}

  Theorem \ref{main_th} implies that  $L(N) \ll 1/ (\log \log \log N)^{C_1}$,
where
${C_1}>0$ is an absolute constant.

  In the present paper we apply  Theorem \ref{main_th} to the theory of dynamical systems
  and obtain result about the multiple recurrence of almost every point
  in an arbitrary dynamical system with two {\it commutative} operators.
  More precisely, we obtain quantitative version of Theorem \ref{Furst1} for
  the case
  $l=2$.

\refstepcounter{section}

{\bf \arabic{section}. On numerical recurrence.}

Let $X$ be a metric space with metric $d(\cdot ,\cdot )$ and a
Borel sigma--algebra of measurable sets $\Phi$. Let $T$ be a
measure preserving transformation of a measure space
$(X,\Phi,\mu)$ and let us assume that measure of $X$ is equal to
$1$. The well--known Poincare (see \cite{Poin}) theorem asserts
that for almost every point $x \in X $:
$$ \forall \eps >  0 ~\forall K  >  0 ~\exists
t  >  K : d(T^tx,x) < \eps .$$
Consider a measure $H_h (\cdot)$ on
$X$, defined as
$$
  H_h (E)= \lim_{\delta \rightarrow 0} H_{h}^{\delta }(E),
$$
where $h(t)$ is a positive ($ h(0)=0 $) continuous increasing
function and $H_{h}^{\delta }(E)= \inf_\tau \{ \sum h(\delta_j)
\}$, when $\tau$ runs through all countable coverings $E$ by open
sets
$ \{ B_j \}$ , $diam(B_i) = \delta_j < \delta$.\\
If $h(t) = t^\alpha$ then we get the ordinary Hausdorff measure
$H_{\alpha } (\cdot ) $.

We shall say that a measure $\mu$ is congruent to a measure $H_h$, if
any $\mu$--measurable set is $H_h$--measurable.

The following theorems
\arabic{section}.1 and \arabic{section}.4
were proved in
\cite{Sh} (see also \cite{Bo,Mo}).

\Th
{\it Let {\it X } be a metric space with
$ H_h (X) = C < \infty $
and let {\it T } be a measure preserving transformation of {\it X}.
Assume that $\mu$ is congruent to $H_h$.
\\Consider the following function:
$
 C(x) = \liminf_{n \to \infty}  \{ n \cdot h(d(T^{n}x,x)) \}
$.
\\
Then the function $C(x)$ is $\mu$--integrable
and for any $\mu$--measurable set $A$
$$
  \int_A C(x) d\mu \le H_h (A).
$$
If $H_h(A)=0$ then $\int_A C(x) d\mu = 0$
with no demand on measures $\mu$ and $H_h$ to be congruent.
}
\label{x_1}

Now we introduce the following concept (see \cite{Kolm}).

\Def
 Let $G$ be a totally bounded subset of $X$.
 By $N_{\eps}(G,X)$ denote the minimal cardinality of $\eps$--net of G.
The number
$H_{\varepsilon }(G,X)$ is called the
{\it $\eps$--entropy} of  $G$.
Put $N_\eps (X) = N_\eps (X,X)$.
\\
If  $X$ is totally bounded then for any $\delta$, we have
$N_\delta (X)< \infty$ and $\sum h(\delta_j) \le N_{\delta
}(X)h(\delta )$. Let $h$ be the function from the definition of
$H_h$. If $ N_{\delta }(X) \le C/h(\delta )$ then $H_{h}(X) \le C
$.

\Def Let $N$ be a natural number.
By $C_{N} (x)$ denote the function
$C_{N}(x) = \min \{~ d(T^{n}x,x) ~|~ 1 \le n \le N ~\}$.
The function $C_{N} (x)$ will be called $N$--constant
of recurrence for point $x$. \label{const_rec}

\Th {\it Let $X$ be  totally bounded metric space with metric
$d(\cdot ,\cdot )$ and function $N(x) = N_x (X)$.
Let $diam (X) = 1$ and ${\it T}$ be a measure--preserving
transformation of {\it X}.
\\
Let $A \subseteq X$ be an arbitrary $\mu$--measurable set and
let $g(x)$ be a real nondecreasing
function bounded on $[0,1]$ such that
for any $t \in (0,1]$ there exists Stieltjes
integral $\int_t^{1} N_{A}(x) dg(x)$,
where $N_{A}(x) = \min (\mu (A),N_x(A,X)/N)$.
Then
$$
  \int_{A} g(C_{N}(x)) d\mu \le
  \inf_t \{ g(t)\mu(A) +
  \int_t^1 N_{A}(x) dg(x) ~\}.
$$
} \label{x_2}

The following lemma due to Poincare (see \cite{Poin,Bo}).
\\
\Lemma {\it Let $Y$ be  $\mu$--measurable set and  $ t \ge 1 $.
Define
$$
  Y(t) := \{ x \in Y ~|~ T^{i}x \notin Y
  \mbox{ for all natural } i, 1 \le i \le t \} .
$$
Then $ \mu(Y(t)) \le 1/t .$
}\label{l_add}

  This lemma is the main tool of the prove of Theorems \ref{x_1}, \ref{x_2}.

  Let us now consider the case of two commutative operators.
  Let
$S$ and $R$ be two {\it commutative} measure--preserving transformation of $X$.
The next result is the main one of this section.

\Th
{\it Let {\it X } be a metric space with
$ H_h (X) = C < \infty $
and let  $S,R$
be two  commutative measure--preserving transformation of $X$.
Assume that $\mu$ is congruent to $H_h$.
\\
Let us consider the function
$$
 C_{S,R}(x) =
 \liminf_{n \to \infty}  \{L^{-1}(n) \cdot max \{ h(d(S^{n}x,x)), h(d(R^{n}x,x)) \} \},
$$
where $L^{-1}(n) = 1/ L(n)$.
\\
Then the function $C_{S,R}(x)$ is $\mu$--integrable
and for any $\mu$--measurable set $A$
$$
  \int_A C_{S,R}(x) d\mu \le H_h (A).
$$
If $H_h(A)=0$ then $\int_A C_{S,R}(x) d\mu = 0$
with no demand on measures $\mu$ and $H_h$ to be congruent. } \label{x_3}

The next definition is analog of Definition \ref{const_rec}.

\Def Let $N$ be a natural number.
By $C^{S,R}_{N} (x)$ denote the function
$C^{S,R}_{N}(x) = \min \{~ max \{ d(S^{n}x,x), d(R^{n}x,x) \} ~|~ 1 \le n \le N ~\}$.
The function $C^{S,R}_{N} (x)$ will be called $N$--constant
of simultaneously recurrence for point $x$.

\Th {\it Let $X$ be a totally bounded metric space
with metric
$d(\cdot ,\cdot )$ and function $N(x) = N_x (X)$.
Let
$diam (X)= 1$ and let ${\it S,R}$ be two measure--preserving transformation of {\it X}.
\\
Let $A \subseteq X$ be an arbitrary  $\mu$--measurable set and
let $g(x)$
be a real nondecreasing
function bounded on $[0,1]$ such that
for any $t \in (0,1]$ there exists Stieltjes
integral $\int_t^{1} N_{A}(x) dg(x)$,
where
 $N_{A}(x) = \min (\mu (A),N_x(A,X) L(N))$.
Then
$$
  \int_{A} g(C^{S,R}_N(x)) d\mu \le
  \inf_t \{ g(t)\mu(A) +
  \int_t^1 N_{A}(x) dg(x) ~\}.
$$
}
\label{x_4}

%
%
%
%

  To prove Theorems \ref{x_3} and \ref{x_4},
  we need several lemmas.

\Lemma  \label{l_ver} {\it Let $ {\bf M} = \{ M_1,\dots ,M_n \}$
be an arbitrary family of
$\mu$--measurable sets.
Let us assume that for any $x \in X $
there exist at most
$l$ sets of the family $\bf{M}$ contain $x$.
Then
$$
  \mu (\bigcup_{i=1}^{n} M_i) \ge \frac{1}{l} \sum_{i=1}^{n} \mu M_i .
$$
}
\Proof
The proof is trivial.

The next lemma is the main of this section.
Using this lemma we obtain Theorems
\ref{x_3} and \ref{x_4} by the same argument as Lemma \ref{l_add} implies Theorems
\ref{x_1} and \ref{x_2}
(for details see \cite{Sh}).

\Lemma {\it Let $Y$ be a $\mu$--measurable set, $ t \ge 1 $.
Define
$$
  Y(t) := \{ x \in Y ~|~ S^{i}x \notin Y \mbox{ or } R^{i}x \notin Y
  \mbox{ for all natural } i, 1 \le i \le t \} .
$$\label{l_key}
Then
$\mu(Y(t)) \le L(t)$.
}
\label{ll}
\\
\Proof
Let $t\ge 1$.
We may assume for convenience that $t$ to be natural.
Define
$$
  M_{k_1, k_2} = S^{-k_1} R^{-k_2} ( Y(t) ) ,~~~~ 1 \le k_1,~k_2 \le t.
$$
Let $x\in X$. By  $A(x)$ denote the set
of indexes
$(k_1, k_2)$ such that $x\in M_{k_1,k_2}$.
Then $A(x) \subseteq [1,t]^2$.
If $|A(x)| > t^2 L(t)$, then using Theorem \ref{main_th},
we obtain that $A(x)$ contains a corner.
Hence there exist $u_1, u_2, u_3 \in Y(t)$
and natural numbers $k,m,d$ such that
$x = S^{-k} R^{-m} u_1 = S^{-k-d} R^{-m} u_2 = S^{-k} R^{-k-d}
u_3$. Since $S$ and $R$ are commutative, it follows that
$u_1 = S^{-d} u_2 = R^{-d} u_3 $.
Hence for $u_1 \in Y(t)$ we
have $S^{d} u_1 \in Y(t)$ and $R^{d} u_1 \in Y(t)$, $d \le t$.
This contradicts the definition  of the set $Y(t)$.
Hence there exist at most $t^2 L(t)$
the  sets  $M_{k_1,k_2}$ such that
$x \in M_{k_1,k_2}$.
Using Lemma \ref{l_ver}, we get
\begin{equation}\label{}
  1 \ge \mu ( \bigcup_{k_1,k_2} M_{k_1, k_2} )
  \ge \frac{1}{t^2 L(t)} \sum_{k_1,k_2} \mu ( M_{k_1, k_2} )
  = \frac{\mu (Y(t))}{L(t)}
\end{equation}
Whence $\mu Y(t) \le L(t)$ as required.

 Now we apply Theorem \ref{x_3} to the case of compact metric space.

The following lemma can be found in \cite{Bg}.

\Lemma {\it Let $X$ be a compact metric space and let
$T_1, \dots T_l$ be
continuous commutative transformations of $X$.
Then there exists a finite measure $\mu$
such that transformations $T_1, \dots T_l$ preserve $\mu$.}

\Cor {\it Let $X$  be a compact metric space with metric $d(\cdot,\cdot)$
and
$H_h (X)< \infty$.
Let $~S, R$ be two continuous commutative transformations of $X$.
Then there exists  $x\in X$ such that
\\
$
 \liminf_{n \to \infty}  \{L^{-1}(n) \cdot max \{ h(d(S^{n}x,x)), h(d(R^{n}x,x)) \} \}
 \le C.
$
}

The author is grateful to Professor N.G. Moshchevitin
  for constant attention to this work.

\end{document}